\documentclass[12pt]{article}

\usepackage{amssymb}
\usepackage{amsmath,amsthm}
\usepackage[dvips]{graphicx}
\usepackage{wrapfig}
\usepackage{url}

\newtheorem{theorem}{Theorem}[section]
\newtheorem{corollary}[theorem]{Corollary}
\newtheorem{lemma}[theorem]{Lemma}
\newtheorem{proposition}[theorem]{Proposition}

\theoremstyle{definition}

\theoremstyle{remark}

\title{The half-twisted splice operation on reduced knot projections}
\author{ Noboru Ito \thanks{Department of Mathematics, Waseda University, 3-4-1 Okubo, Shinjuku-ku, Tokyo 169-8555, Japan. noboru@moegi.waseda.jp}
\and Ayaka Shimizu \thanks{Department of Mathematics, Hiroshima University, 1-3-1 Kagamiyama, Higashi-Hiroshima, 739-8526, Japan. shimizu1984@gmail.com}}
\date{\today}

\begin{document}

\maketitle

\begin{abstract}
We show that any nontrivial reduced knot projection can be obtained from a trefoil projection 
by a finite sequence of half-twisted splice operations and their inverses such that the result of each step in the sequence is reduced. 
\end{abstract}

\section{Introduction}

Throughout this paper, knot projections are on $S^2$. 
A \textit{knot projection} is a regular projection of a knot. 
Arnold \cite{arnold}, Calvo \cite{calvo}, and Endo-Itoh-Taniyama \cite{EIT} introduced a local move that they, respectively, called perestroika, twisted splice, and smoothing operation, to describe a connection between plane curves or knot projections or diagrams. 
(See also \cite{arnold-e}, \cite{arnold2}, \cite{calvo2}, \cite{hoste}, \cite{taniyama}, \cite{taniyama2}.) 
In this paper, we study a half-twisted splice on knot projections to create a connection between \textit{reduced} knot projections. \\

\begin{figure}[h]
 \begin{center}
  \includegraphics[width=47mm]{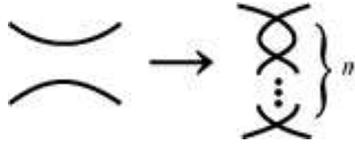}
 \end{center}
 \caption{Twisted splice.}
 \label{twisted-s}
\end{figure}

\begin{figure}[th]
 \begin{center}
  \includegraphics[width=125mm]{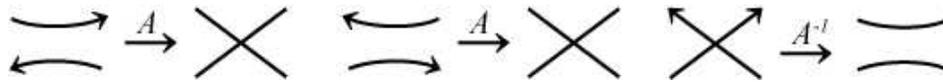}
 \end{center}
 \caption{Half-twisted splice.}
 \label{def-a}
\end{figure}

A \textit{twisted splice} on a knot projection is the replacement of trivial tangle with a two-arc braided tangle with $m$ crossings resulting in a knot projection as shown in Fig. \ref{twisted-s}, where $m$ is even (resp. odd) for segments oriented in the same (resp. opposite) direction \cite{calvo}. 
We call the twisted splice with $m=1$, shown in Fig. \ref{def-a}, a \textit{half-twisted splice}, and denote it by $A$. 
Obviously $A$ is a local move on knot projections. 
The inverse $A^{-1}$ is shown in the right-hand side of Fig. \ref{def-a}. 
Note that a twisted splice does not depend on the orientation of knot projections, but 
depends only on the local ``relative orientation''. 
Hence $A$ and $A^{-1}$ also depend only on the relative orientation. \\

\begin{figure}[th]
 \begin{center}
  \includegraphics[width=117mm]{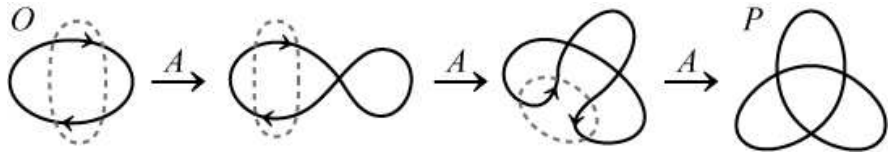}
 \end{center}
 \caption{Half-twisted splices.}
 \label{tre-a}
\end{figure}

The \textit{trivial projection} is the knot projection $O$, with no crossings, as depicted in the left-hand side of Fig. \ref{tre-a}. 
The \textit{trefoil projection} is the knot projection $P$ as depicted in the right-hand side of Fig. \ref{tre-a}. 
We can obtain the trefoil projection by three half-twisted splices from the trivial projection as shown in Fig. \ref{tre-a}. 
Every knot projection $P$ with $n$ crossings can be obtained from a trivial projection by $n$ half-twisted splices because we obtain a knot projection with $(n-1)$ crossings from $P$ by a single $A^{-1}$. 
(This is suggested by Kawauchi, and see also \cite{calvo}.) \\

\begin{figure}[h]
 \begin{center}
  \includegraphics[width=40mm]{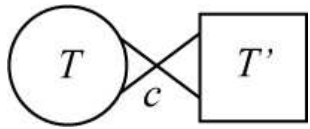}
 \end{center}
 \caption{Reducible projection.}
 \label{red}
\end{figure}

A knot projection $P$ is \textit{reducible} and has a \textit{reducible crossing} $c$ if $P$ has a crossing point $c$ as shown in Fig. \ref{red}, where $T$ and $T'$ imply tangles. 
For example, the knot projections in the middle in Fig. \ref{tre-a} are reducible projections. 
A knot projection $P$ is \textit{reduced} if $P$ is not reducible. 
Note that around a reducible (resp. non-reducible) crossing, there are exactly three (resp. four) disjoint regions. \\

In this paper, we shall prove the following theorem: 

\phantom{x}
\begin{theorem}
Every nontrivial reduced knot projection is obtained from the trefoil projection by a finite sequence of $A$ and $A^{-1}$ such that the result of each step in the sequence is reduced. 
\label{mainthm1}
\end{theorem}
\phantom{x}

\noindent This immediately gives us the following corollary: 

\phantom{x}
\begin{corollary}
Let $P$, $Q$ be nontrivial reduced knot projections. 
There is a finite sequence of $A$ and $A^{-1}$ from $P$ to $Q$ such that the result of each step in the sequence is reduced. 
\end{corollary}
\phantom{x}

This study was motivated by the knot game Region Select which was created by Kawauchi, Kishimoto and Shimizu. 
(This game is based on a local move, region crossing change \cite{shimizu}.) 
In the game, we use knot projections. 
In particular, reduced knot projections make the game more interesting. 
A preliminary web-based version of Region Select in 2011 used $A$ to obtain a reduced knot projection from a reduced one. 
(The present version of the game can be found at \\
\url{http://www.sci.osaka-cu.ac.jp/math/OCAMI/news/gamehp/etop/gametop.html}.)\\

The rest of this paper is organized as follows: 
In Section \ref{sec-anti-p}, we introduce two new moves $B$ and $C$ and explore their relationship to $A$ and $A^{-1}$. 
In Section \ref{sec-pf-thm1}, we prove Theorem \ref{mainthm1} by looking at the set of Reidemeister moves appropriate to oriented knot projections.

\section{The half-twisted splice and moves $B$ and $C$}
\label{sec-anti-p}

\begin{figure}[h]
 \begin{center}
  \includegraphics[width=30mm]{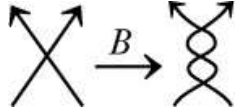}
 \end{center}
 \caption{Definition of move $B$.}
 \label{t23}
\end{figure}

In this section, we discuss what we can accomplish with the operations $A$ and $A^{-1}$. 
We define the move $B$ to be the local move of a knot projection at a crossing which adds exactly two crossings as depicted in Fig. \ref{t23}. \\

We have the following proposition: 

\phantom{x}
\begin{proposition}
Let $P$ be a knot projection, and $c$ a crossing of $P$. 
Let $P'$ be the knot projection obtained from $P$ by $B$ at $c$. 
We can obtain $P'$ from $P$ by a finite sequence of $A$ and $A^{-1}$. 
Further, if $P$ is reduced, the result of each step in this sequence is also reduced. 

\begin{proof}
See Fig. \ref{t23-aa}. 
If $P$ is reduced, then the regions $R_1$, $R_2$, $R_3$ and $R_4$ around $c$ are pairwise disjoint. 
Hence each step in the sequence of the local moves never create a reducible crossing. 
\begin{figure}[th]
 \begin{center}
  \includegraphics[width=130mm]{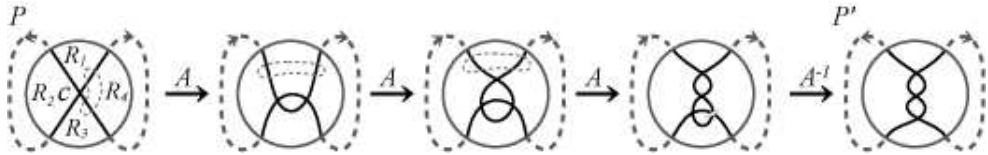}
 \end{center}
 \caption{$B$ is realized by $A$ and $A^{-1}$.}
 \label{t23-aa}
\end{figure}
\end{proof}
\label{t23-aa-p}
\end{proposition}
\phantom{x}

\begin{figure}[th]
 \begin{center}
  \includegraphics[width=80mm]{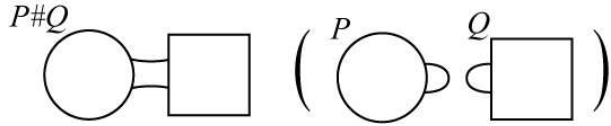}
 \end{center}
 \caption{Connected sum.}
 \label{conn}
\end{figure}

\begin{figure}[th]
 \begin{center}
  \includegraphics[width=75mm]{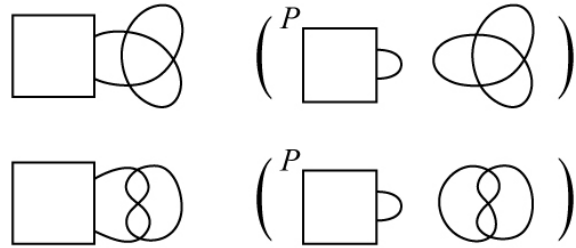}
 \end{center}
 \caption{Connected sums with trefoil projection.}
 \label{two-conn}
\end{figure}

Next we consider connected sums. 
A \textit{connected sum} $P\sharp Q$ of two knot projections $P$ and $Q$ is a knot projection as depicted in Fig. \ref{conn}. 
Note that the connected sum is not unique. 
For example, the two connected sums in Fig. \ref{two-conn} are connected sums of a knot projection $P$ and a trefoil projection. 
The upper (resp. bottom) one is obtained by the connected sum at a triangle (resp. bigon) region of the trefoil projection. 
We have the following proposition: 

\begin{proposition}
The two types of connected sums with a trefoil projection are equivalent by a finite sequence of $A$ and $A^{-1}$ which create no reducible crossing. 

\phantom{x}
\begin{proof}
See Fig. \ref{conn-pf2}. 
Each crossing in Fig. \ref{conn-pf2} is non-reducible, and each box represents the same tangle. 
\begin{figure}[th]
 \begin{center}
  \includegraphics[width=100mm]{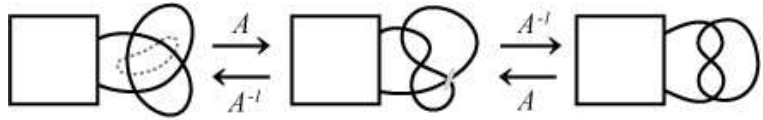}
 \end{center}
 \caption{Connected sums with a trefoil projection.}
 \label{conn-pf2}
\end{figure}
\end{proof}
\end{proposition}
\phantom{x}

\begin{figure}[th]
 \begin{center}
  \includegraphics[width=35mm]{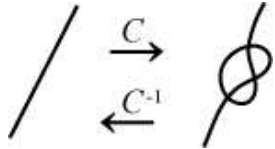}
 \end{center}
 \caption{Definition of move $C$.}
 \label{f-move}
\end{figure}

\noindent Now we define moves $C$ and $C^{-1}$ to be the local moves shown in Fig. \ref{f-move}. 
Namely, $C$ is the connected sum with a trefoil projection, and $C^{-1}$ is the ``removing'' a trefoil projection. \\

We have the following proposition:

\phantom{x}
\begin{proposition}
We can create $C$ anywhere on a non-trivial knot projection by a finite sequence of $A$ and $A^{-1}$. 
Further, if $P$ is reduced, the result of each step in this sequence is also reduced. 

\begin{proof}
See Fig. \ref{tre-m}. 
If $P$ is reduced, then the regions $R_1$, $R_2$, $R_3$ and $R_4$ are pairwise disjoint, and therefore all the crossings shown in Fig. \ref{tre-m} are non-reducible. 
\begin{figure}[th]
 \begin{center}
  \includegraphics[width=120mm]{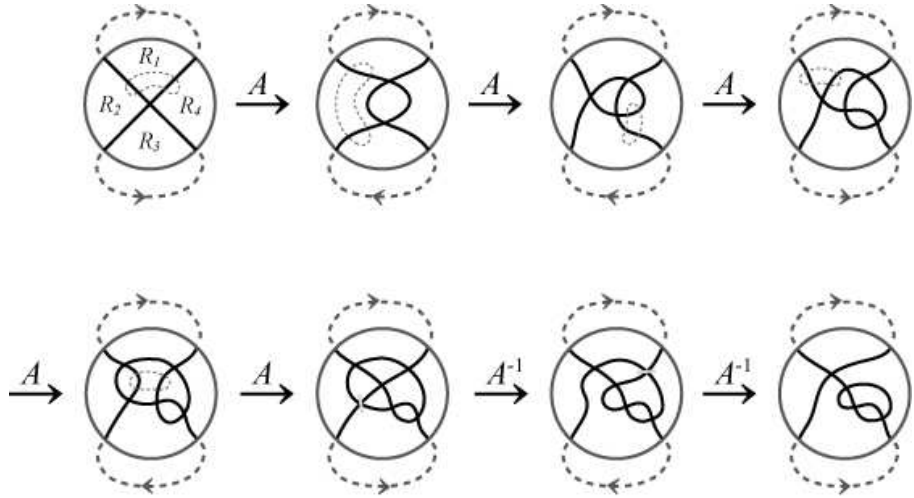}
 \end{center}
 \caption{$C$ is realized by $A$ and $A^{-1}$.}
 \label{tre-m}
\end{figure}
\end{proof}
\label{aa-f}
\end{proposition}
\phantom{x}

\noindent By repeated application of Proposition \ref{aa-f}, we have the following corollary: 

\phantom{x}
\begin{corollary}
Any connected sum of $n$ trefoil projections $(n\ge 2)$ is obtained from a trefoil projection 
by a finite sequence of $A$ and $A^{-1}$ such that the result of each step in the sequence is reduced. 
\label{con-t-aa}
\end{corollary}
\phantom{x}

\section{Proof of Theorem \ref{mainthm1}}
\label{sec-pf-thm1}

\begin{figure}[th]
 \begin{center}
  \includegraphics[width=100mm]{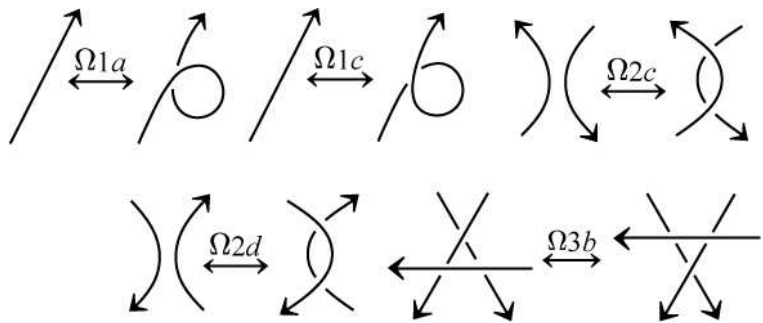}
 \end{center}
 \caption{Reidemeister moves.}
 \label{ori-rei}
\end{figure}

In this section, we prove Theorem \ref{mainthm1}. 
First, we consider Reidemeister moves on oriented knot diagrams and projections by refering to Polyak's theorem. 
Polyak \cite{polyak} showed that, if $D$ and $D'$ are diagrams in $\mathbb{R}^2$ representing the same oriented link, then $D'$ is obtained from $D$ by a finite sequence of oriented Reidemeister moves $\Omega _{1a}$, $\Omega _{1c}$, $\Omega _{2c}$, $\Omega _{2d}$, and $\Omega _{3b}$ as shown in Fig. \ref{ori-rei}. \\

\noindent By erasing the over- and under-crossings in the diagram, we obtain ten distinct Reidemeister moves on oriented projections. 
Remark that every knot projection can be made the projection of an unknot by an appropriate choice of crossings. 
In this setting, Polyak's theorem becomes: 

\phantom{x}
\begin{corollary}
Any oriented knot projection can be obtained from the trivial projection by a finite sequence of the local moves 
$\Omega _{1+}$, $\Omega _{1-}$, $\Omega _{2+}$, $\Omega _{2-}$, $\Omega _{2'+}$, $\Omega _{2'-}$, $\Omega _{3+}$ and $\Omega _{3-}$ shown in Fig. \ref{ori-rei-ori}. 
\begin{figure}[th]
 \begin{center}
  \includegraphics[width=125mm]{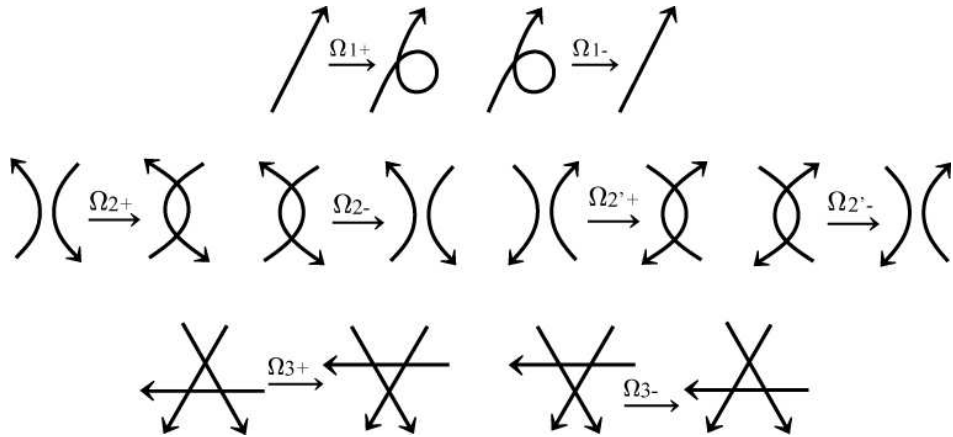}
 \end{center}
 \caption{Reidemeister moves on projections.}
 \label{ori-rei-ori}
\end{figure}
\label{ori-123}
\end{corollary}
\phantom{x}

\noindent Let $P$ be a knot projection. 
Let $P_B$ be the knot projection obtained from $P$ by applying $B$ at all the crossings of $P$. 
We have the following proposition: 

\phantom{x}
\begin{proposition}
$P_B$ is a reduced knot projection even if $P$ is reducible. 

\phantom{x}
\begin{proof}
Each reducible crossing of $P$ will be replaced with the three non-reducible crossings as shown in Fig. \ref{bp-red}. 
\begin{figure}[th]
 \begin{center}
  \includegraphics[width=75mm]{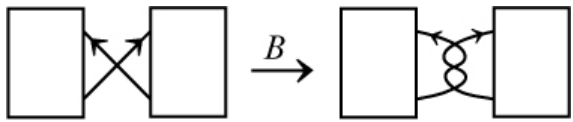}
 \end{center}
 \caption{A reducible crossing will be replaced with non-reducible crossings.}
 \label{bp-red}
\end{figure}
\end{proof}
\end{proposition}

\noindent We have the following lemma: 

\phantom{x}
\begin{lemma}
$P_B$ is obtained from the trefoil projection by a finite sequence of $A$ and $A^{-1}$ without ever passing through a reducible projection. 

\phantom{x}
\begin{proof}
The moves $B$ and $C$ and their inverses are allowed in this proof since they each consist of a sequence of $A$ and $A^{-1}$. 
By Corollary \ref{ori-123}, there exists a finite sequence of knot projections 
$P=P^0\rightarrow P^1\rightarrow P^2\rightarrow \dots \rightarrow P^m=O$, 
where each $\rightarrow$ is one of 
$\Omega _{1+}$, $\Omega _{1-}$, $\Omega _{2+}$, $\Omega _{2-}$, $\Omega _{2'+}$, $\Omega _{2'-}$, $\Omega _{3+}$ and $\Omega _{3-}$. 
Remark that the last step of the sequence is $\Omega _{1-}$, $\Omega _{2-}$ or $\Omega _{2'-}$. 
We can replace the $\Omega _{2-}$ move with a pair of $\Omega _{1-}$ moves if the last step is $\Omega _{2-}$ move. 
Similarly, we can replace the $\Omega _{2'-}$ move with a pair of $\Omega _{1-}$ moves by reversing the orientation of the knot projection if the last step is $\Omega _{2'-}$ move. 
Hence we can assume that $P^{m-1}$ is the knot projection with exactly one crossing. 
Let $P_B^i$ be the knot projection obtained from $P^i$ by applying $B$ at all the crossings of $P^i$ ($i=0,1,\dots ,m$). 
Note that $P^{m-1}_B$ is the trefoil projection. 
Now we show that we can move from $P_B^i$ to $P_B^{i+1}$ by a finite sequence of $A$ and $A^{-1}$ without ever passing through a reducible projection ($i=0,1,\dots ,m-1$).  \\

\begin{description}
 \item[\bf Case 1: Reidemeister move $\Omega _{1+}$. ] 
Use move $C$ on $P_B^{i}$ as in Fig. \ref{case-one}. 
\begin{figure}[h]
 \begin{center}
  \includegraphics[width=60mm]{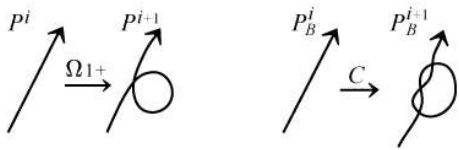}
 \end{center}
 \caption{$\Omega _{1+}$ with $B$ is realized by $C$.}
 \label{case-one}
\end{figure}

\item[\bf Case 2: Reidemeister move $\Omega _{1-}$. ]
Use move $C^{-1}$ on $P_B^{i+1}$ as in Fig. \ref{case-two}.
\begin{figure}[h]
 \begin{center}
  \includegraphics[width=60mm]{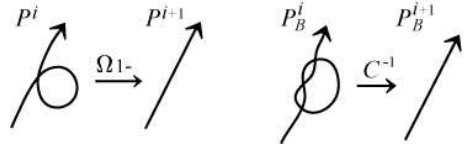}
 \end{center}
 \caption{$\Omega _{1-}$ with $B$ is realized by $C^{-1}$.}
 \label{case-two}
\end{figure}

\item[\bf Case 3: Reidemeister move $\Omega _{2+}$. ] 
See Fig. \ref{sig-2}. 
Remark that $P^i_B$ is a nontrivial projection. 
Hence we can apply $C$ to the projection. 
Note that all the crossings shown in the figure are non-reducible crossings. 
\begin{figure}[h]
 \begin{center}
  \includegraphics[width=120mm]{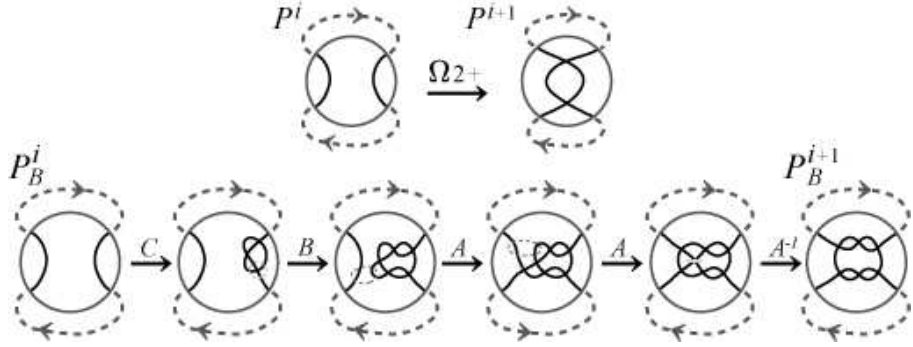}
 \end{center}
 \caption{$\Omega _{2+}$ with $B$ is realized by $A$ and $A^{-1}$.}
 \label{sig-2}
\end{figure}

\item[\bf Case 4: Reidemeister move $\Omega _{2-}$. ]
See Fig. \ref{sig-2} in the reverse order.

\item[\bf Case 5: Reidemeister move $\Omega _{2'+}$] 
Consider $\Omega _{2+}$ with the orientation of a projection reversed. 

\item[\bf Case 6: Reidemeister move $\Omega _{2'-}$. ]
Similarly, consider $\Omega _{2-}$ with the orientation reversed. 

\item[\bf Case 7: Reidemeister move $\Omega _{3+}$. ]
See Fig. \ref{type-a} and \ref{type-b}. 
Note that all the crossings shown in the figures are non-reducible crossings. 
\begin{figure}[hp]
 \begin{center}
  \includegraphics[width=120mm]{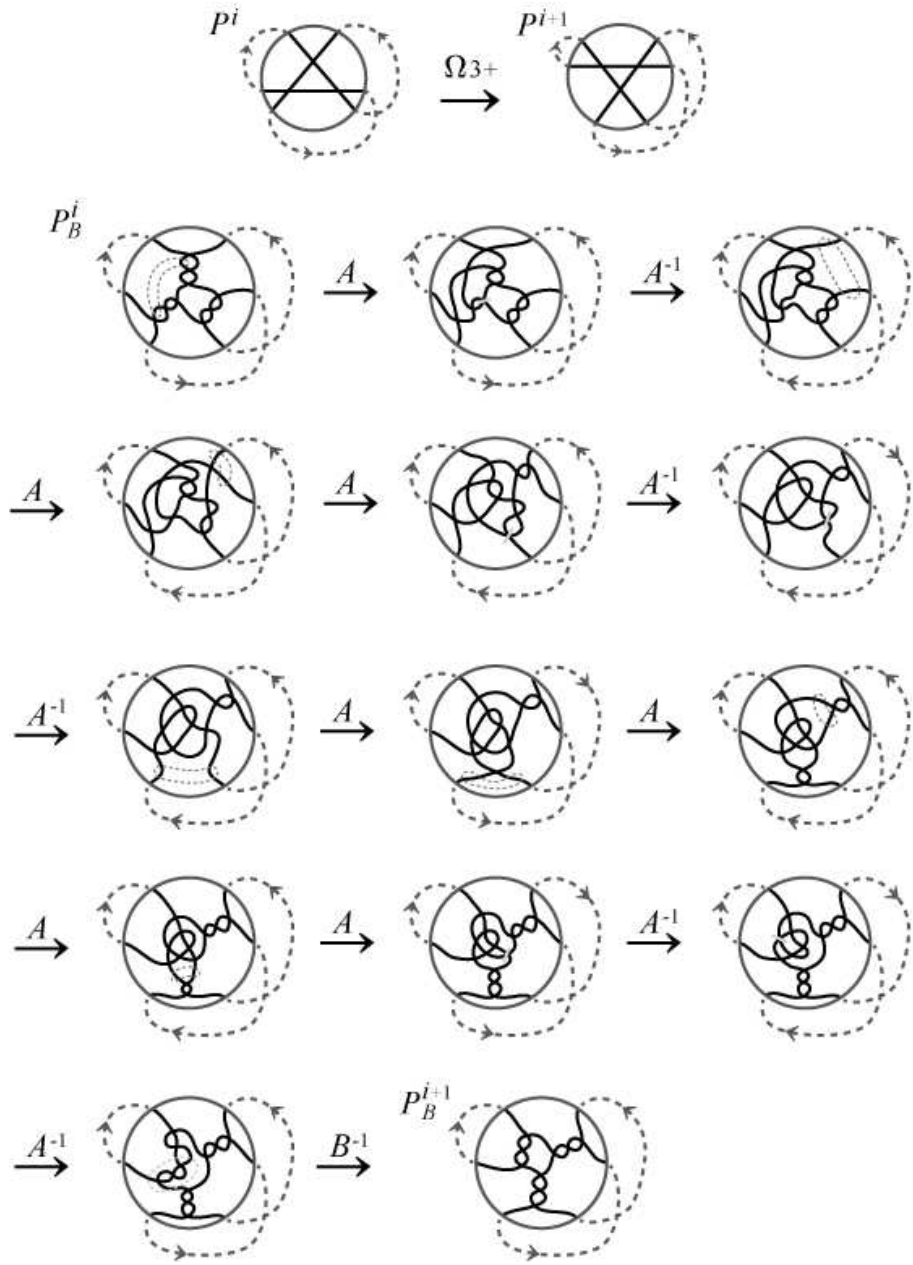}
 \end{center}
 \caption{$\Omega _{3+}$ with $B$ is realized by $A$ and $A^{-1}$ I.}
 \label{type-a}
\end{figure}
\begin{figure}[hp]
 \begin{center}
  \includegraphics[width=120mm]{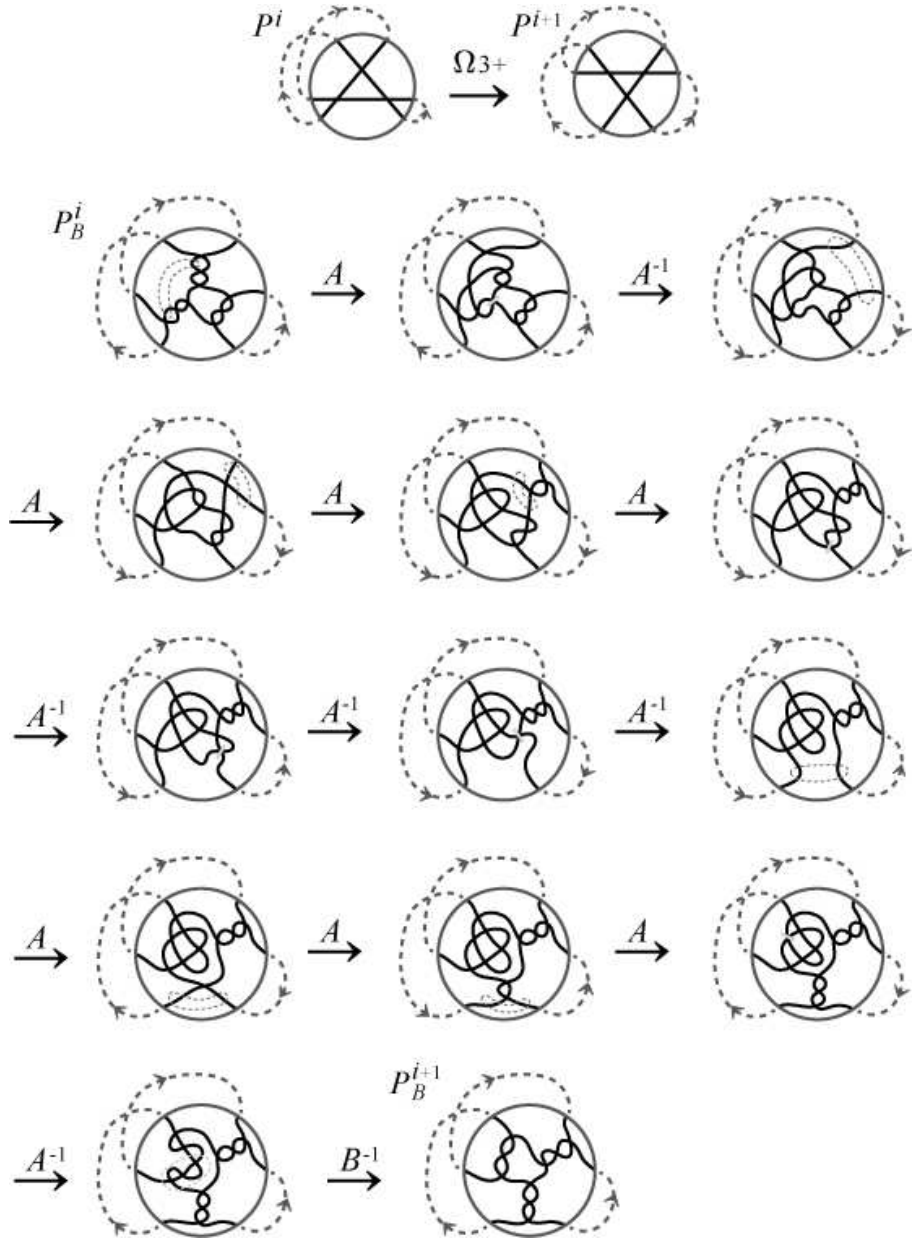}
 \end{center}
 \caption{$\Omega _{3+}$ with $B$ is realized by $A$ and $A^{-1}$ II.}
 \label{type-b}
\end{figure}

\item[\bf Case 8: Reidemeister move $\Omega _{3-}$. ]
Similar to the case $\Omega _{3+}$. 
\end{description}
\end{proof}
\label{aaf-con}
\end{lemma}
\phantom{x}

\noindent We prove Theorem \ref{mainthm1}: 

\phantom{x}
\noindent \textit{Proof of Theorem \ref{mainthm1}.} 
We can transform $P$ into $P_B$ by Proposition \ref{t23-aa-p}. 
Then, we can transform $P_B$ by a finite sequence of $A$ and $A^{-1}$ into the trefoil projection by Lemma \ref{aaf-con}. 
In each case, the result of each step in the sequence will be a reduced projection. \\
\hfill$\square$\\
\phantom{x}

\section*{Acknowledgment}
N.I. was supported by Grant-in-Aid for Young Scientists (B) (23740062). 
A.S. thanks the members of Friday Seminar on Knot Theory 2011 at Osaka City University for valuable discussions and advice. 
In particular, she thanks Professor Akio Kawauchi for many helpful suggestions and giving her information about local moves on knot projections. 
She is deeply grateful to Professor Kenneth C. Millett for valuable discussions, advice, encouragement, and introducing her to Calvo's paper. 
She also thanks the referee for the kind advice and suggestions. 
She was partly supported by JSPS Research Fellowship for Young Scientists.

\end{document}